\documentclass[12pt]{article}
\pdfoutput=1

\usepackage{graphicx}
\usepackage{color}
\usepackage{amsfonts}
\usepackage{amssymb}
\usepackage{amsmath}
\usepackage{fullpage}

\def\rstyle#1{\mbox{\small #1}}

\def\avi{{\rstyle{av}(i)}}
\def\ad{\rstyle{ad}}
\def\wat{\rstyle{wat}}
\def\gen{\rstyle{gen}}
\def\tt#1#2{t_{#1\leftarrow #2}}

\def\Vmin#1{V_{#1}^{\rstyle{min}}}
\def\Vmax#1{V_{#1}^{\rstyle{max}}}
\def\Tmin#1{T_{#1}^{\rstyle{min}}}
\def\Tmax#1{T_{#1}^{\rstyle{max}}}
\def\Dmax#1{D_{#1}^{\rstyle{max}}}
\def\eqdef{\stackrel{\rstyle{def}}{=}}
\def\uad{{\cal T}^{\rstyle{ad}}}
\def\U2ad{{\cal U}^{\rstyle{ad}}}

\def\inutile#1{}

\def\CfinI{C^{\rstyle{end}}_i}
\def\psca#1{{\langle  #1 \rangle}}

\def\ukp1{u^{k+1}}

\def\ukp1{u^{k+1}}

\def\jversi{{i\leftarrow j}}

\def\T#1#2{{\Theta}_{#1}(#2)}

\def\norm#1{\Vert #1 \Vert }
\let\dis=\displaystyle
\def\ds#1{{\displaystyle #1}}
\def\ligne #1 & #2 & #3{{\displaystyle  #1} & #2 & {\displaystyle #3} }
\def\text#1{\quad\mbox{#1}\quad}

\def\dynamiqueu{V_i(t+1) = V_i(t) - D_i(t) - T_i(t) + A^{*}_i(t)}
\def\dynamiqued{\sum_{j \in \Gamma^T_i} \left(D_j(t_\jversi)
			+T_j(t_\jversi)\right)}
\def\dynamique{\hspace{- 0.4cm} \dynamiqueu \nonumber \\ &&\hspace{-0.1cm}+ \dynamiqued}

\def\touttd{\forall t \in [0,H[}
\def\touti{\forall i \in {\cal I}}
\def\exclusion#1#2{\left(#1_i(t) -\Vmax{i} \right) #2_i(t) =0 }

\def\turbi#1{#1_i([0,H[) \in \uad_i}
\def\critere{\sum_{t \in [0,H[} \hspace{-0.2cm} - c_i(T_i(t),t) - \CfinI(V_i(H),H)}
\def\criteredans{- c_i(T_i(t),t) - \CfinI(V_i(H),H)}
\def\sument{\sum_{t \in [0,H[}}
\def\gain{\sum_{t \in [0,H[}  c_i(T_i(t),t) + \CfinI(V_i(H),H)}
\def\Vbornes#1{#1_i(t) \in [\Vmin{i},\Vmax{i}]}
\def\Tbornes#1{#1_i(t) \in [\Tmin{i},\Tmax{i}]}
\def\Dbornes#1{#1_i(t) \in [0,\Dmax{i}]}
\def\bornes#1#2#3{\begin{array}{c}
\Vbornes{#1}\,,\\*\;\Tbornes{#2}\,, \\*
\Dbornes{#3}\,,
\end{array}}

\def\Bbb#1{\mathbb{#1}}

\date{1994} 
\author{J-P. Chancelier \\  ENPC,  CERGRENE,\\ 
	Noisy le Grand, France \and A.Renaud \\  
	Electricit\'{e} de France,\\  Clamart, France }

\title{Daily Generation Scheduling~:\break Decomposition Methods to
		Solve the Hydraulic Problems}
\begin{document}
\maketitle

\begin{abstract}
  Short-term hydro-generation management poses a non-convex
  or even non-continuous optimization problem.
  For this reason, the problem of systematically obtaining feasible and
  economically satisfying solutions has not yet been completely solved.

  Two decomposition methods, which, as far as we know, have not been
  applied in this field, are here proposed~: 
  \begin{itemize}
  \item the first is based on a decomposition by prediction
    method and  the coordination is a
    primal-dual relaxation algorithm,
  \item handling the dynamic constraints by duality,
    the second achieves a price decomposition by
    an Augmented Lagrangian technique.
  \end{itemize} 
  Numerical tests show the efficiency of these algorithms. They will
  enable the process in use at Electricit\'{e} de France to be improved.
\end{abstract}

\subsection*{Keywords} Hydro-thermal scheduling, Decomposition methods.
\section{Introduction} 
The optimization of a hydro-valley's daily generation schedules poses a problem of
an appreciable size. For instance, a problem related to a valley of
five or six reserves,  with a
half-hour step has no less than 250 or 300 constraints.

Since the 1970's, very efficient linear programming methods and softwares  have been designed to handle such
optimization problems. The most commonly used
algorithms take  advantage of the
network flow structure of the dynamic constraints~\cite{maurras} and
succeed in being about  a hundred
times faster than standard linear programming methods~\cite{ovide}. 

Over the last ten years, extensions to the nonlinear (but convex) case have been made. Using
Frank-Wolfe or projected reduced gradient techniques, efficient softwares
have been developed. Thanks to them, nonlinear efficiency curves, for
example,  can be coped with.

Nevertheless, the modelling of the hydro-problems which can be solved by these methods is not
completely satisfactory. Downstream flow requirements or ``spillage
constraints'' cannot be handled. Moreover, the efficiency curves,
which are generally non-convex, have to be roughly
approximated and no discontinuity in the generating domain can be dealt with. 

\vskip 0.1 cm

The Electricit\'{e} de France generation mix has over $150$ thermal groups and $15$ valleys.
The global optimization of the daily generation schedules is achieved by using a price decomposition
method. In this optimization, some reservoirs are aggregated and the generating domain of the
hydro-plants is assumed to be convex. Considering the number of local hydro-problems which have
to be solved during this optimization, such a simplifying hypothesis can easily be understood.  

Nevertheless, at regional level, a second optimization stage is necessary~\cite{PSCC86}. 
For each valley, a schedule is computed independently~:
\begin{itemize}
\item firstly, a linear problem is solved. It takes into
  account the dual variables that the
  global optimization yields and handles a detailed but convex
  modelling of the constraints,
\item secondly, in the neighborhood of this schedule, a heuristic
  ``smoothing'' software processes a feasible solution with
  respect to the non-convex or even non-continuous constraints.
\end{itemize} 

The study, whose results are here outlined, aims at improving this regional
level two-steps process. Two decomposition strategies have been tested in order
to solve this problem directly (in one stage). These algorithms only require
small-sized nonconvex sub-problems to be solved at each iteration. Therefore, an
exact method, such as dynamic programming, can be used to deal with this local
problems.

\vskip 0.1 cm

The paper is organized as follows. In section 2, we formulate the
considered optimization problem. Then, in section 3, we outline the
theoretical background of the prediction and price decomposition
strategy along with their application to the hydro-problem. Finally,
in section 4, we present the results of the numerical tests which have been realized.
\section{The Hydro Problem}
\label{formulation}
\begin{figure}[hbtp]
  \begin{center}
    \includegraphics{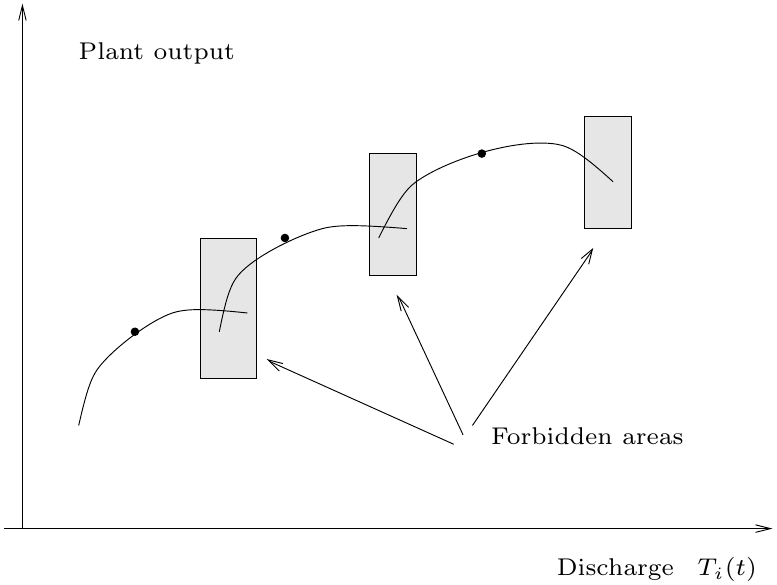}
  \end{center}
  \caption{\label{noncont}The generating domain of the plant $i$ is non-continuous}
\end{figure}

\begin{figure}[hbtp]
  \begin{center}
    \includegraphics{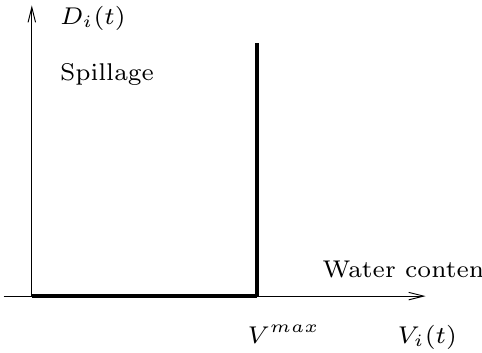}
  \end{center}
  \caption{\label{debord}The domain defined by the constraints $\exclusion{V}{D}$ is nonconvex}
\end{figure}

As already emphasized, the ``regional'' problem we are considering consists in the schedule
optimization of one valley. It may be written as follows~:
\begin{eqnarray} 
  \hspace{-0.2cm} \lefteqn{ \displaystyle \min_{T_i(.)}   \sum_{i \in {\cal I}}
  \critere } \label{pb}\\*
  \lefteqn{ \displaystyle \mbox{subject to~:}}  \nonumber\\*
  \lefteqn{ \displaystyle \bullet \;\touti\,,\;\touttd~:}	 \nonumber \\*
  &&\dynamique \label{dynam}\\*
  &&\exclusion{V}{D} \label{exclus}\\*
  &&\left[ \bornes{V}{T}{D} \right. \label{borne} \\*
  \lefteqn{ \displaystyle \bullet\; \touti~:  \turbi{T} }\label{turbi}
\end{eqnarray}
where~:
\begin{itemize}
\item ${\cal I}$ is the set of water reservoirs of the valley under
  consideration --- each reservoir is related to a plant having the same index ---, 
\item $[0,H]$ is the studied period,
\item $V_{i}(t)$ is the water content of the reservoir $i$ at time $t$ and $\Vmin{i}$, $\Vmax{i}$
  are respectively the upper and the lower bounds of the reservoir $i$,
\item $T_{i}(t)$ is the discharge of plant $i$ over
  $[t,t+1[$ and $\Tmin{i}$ et $\Tmax{i}$ are respectively the lower and the upper
  bounds of the discharge over the period,
\item $\uad_i$ is defined by the generating constraints of the plant $i$, 
\item $D_{i}(t)$ is the spillage of plant $i$
  over $[t,t+1[$ and $\Dmax{i}$ its upper bound over the period,
\item $\Gamma^T_{i}$ is the set of plants located upstream $i$ --- whose discharges are inflows of
  $i$ ---, 
\item $\avi$ is the reservoir located just downstream~$i$,
\item $t_{\jversi} \eqdef t - \delta_{\jversi}$ with
  $\delta_{\jversi}$ denoting the delay
  for discharge of plant~$j$ to reach reservoir~$i$,
\item $A^*_{i}(t)$ is the natural inflows of the reservoir~$i$
  during the period $[t,t+1[$. These inflows
  are supposed to be known, 
\item $\CfinI(V_i(H),H)$ is the water value of the reservoir $i$ at the end of the studied period,
\item $c_{i}(T_i(t),t)$ is the ``value'' of the generation
  related to the discharge $T_i(t)$ at time $t$. 
\end{itemize} 
It has already been pointed out that the 
generating constraints we consider define a non-continuous domain 
$\uad_i$. We illustrate this in Figure~\ref{noncont}. Moreover, as is
shown  in Figure~\ref{debord} the constraints~(\ref{exclus}) also introduce nonconvexities.

The numerical tests, hereafter outlined, use the current EDF regional modelling. Concerning the
generating constraints $\uad_i$~:
\begin{itemize}
\item only a finite number of values is allowed for the discharge $T_i(t)$ (See Figure~\ref{fig:usine}),
\item there is a minimum delay between two variations of the
  discharge. Moreover, these variations have to be smooth.
\end{itemize}
The cost function $c_i(.,t)$ is a piecewise\-linear function of the discharge and 
$\CfinI(V_i(H),H)$
is a linear function of the final content~: 
$$\CfinI(V_i(H),H) = c_{i}^{\wat} \left(V_i(H) - V_{i}(0) \right)\enspace.$$
This modelling allows us to compare decomposition methods to the current EDF process. Nevertheless,
the decomposition methods framework hereafter presented, could indeed
be applied to more general modelling

\section{Resolution Methods}
In order to present the coordination-decomposition methods we use, we will consider the following
problem:
\begin{equation}
  \displaystyle \min_{u \in \U2ad} J(u) \enspace,\enspace \mbox{subject to:} \enspace\Theta(u)=0
  \enspace, \label{pbtyp}
\end{equation}
where $\U2ad$ is a closed set of a vector space ${\cal U}$, $J:{\cal U}  \rightarrow
{\Bbb R}$ and $\Theta: {\cal U} \rightarrow {\cal C}$ an
affine function from ${\cal U}$ to the vectorial space ${\cal C}$. To allow
decomposition, we also assume~:
\begin{eqnarray}
  {\cal U} = \prod_{i \in I} {\cal U}_i \enspace,\enspace
  \U2ad = \prod_{i \in I} \U2ad_i \enspace, \nonumber \\*
  \forall u \in \U2ad\enspace :\enspace J(u)=\sum_{i \in I} J_i(u_i)\enspace, \nonumber
\end{eqnarray}
where $I$ is a finite set, $\U2ad_i \subset {\cal U}_i$ and $ J_i: {\cal U}_i \rightarrow {\Bbb R}  $. 
\subsection{Today's resolution method}
First of all, the convexified problem~(problem~(\ref{pb}), without
constraints (\ref{exclus}) and (\ref{turbi})) is solved by linear
programming. Secondly, following the course of the river,
for every plant, a feasible schedule is processed. In this aim, taking into account {\it all} the 
generating constraints\footnote{The inflows are known~: feasible
  discharges of the upstream plants have  already been processed.
}, discharge is computed by minimizing a mean-square distance to the solution of the
convexified problem. These subproblems are solved by dynamic programming. 

It may be noticed that this heuristic has a major drawback~: it does not ensure that a feasible 
solution will be
found. The inflows being given, a subproblem may have no solution. 

\subsection{Price decomposition} 
\subsubsection{Theoretical background} 
Provided that  the cost function is separable, the Uzawa algorithm \cite{uzawa} is certainly the most commonly
used to achieve this type of decomposition. If we suppose $\Theta(u)=
\sum_{i \in I}\Theta_i(u_i)$,then,  applied to the problem~(\ref{pb}), this algorithm may be outlined as
follows ($k$~iteration): 
\begin{eqnarray}
  \forall i \in I\;:\; \dis \min_{u_i \in \U2ad_i} J_i(u_i) + \psca{ p^k,
  \Theta_i(u_i)} \; \Rightarrow \; u^{k+1}_i  \nonumber \\*
  p^{k+1} = p^k + \varepsilon^k \Theta(u^{k+1}) \hfill  \nonumber
\end{eqnarray}
Therefore, at each iteration the $i$-subsystem minimizes a balance (the Lagrangian) which takes into
account its own cost function $J_i$ and a ``revenue'' $\psca{ p^k,
  \Theta_i(u_i)}$, from its contribution to the satisfaction of the constraint $\Theta(u)=0$.

More formally, the Lagrangian related to (\ref{pbtyp}) is defined over $ {\cal U} \times {\cal C} $ as follows~:
\begin{eqnarray}
  \forall (u,p) \enspace : \enspace 
  L(u,p) \eqdef J(u) +  \psca{ p, \Theta(u)} \enspace. \label{lagrang}
\end{eqnarray}
The above mentioned algorithm may be understood as maximizing the
dual function:
\[
  \Psi: p \mapsto \min_{u \in U^{\ad}} L(u,p) 
\]
by a gradient type algorithm. 

If the function $J$ is not strictly convex, as in the hydro-problem we are dealing with, this dual
function is not differentiable. Consequently, to ensure the convergence, a sub-gradient algorithm
must be used to maximize $\Psi$. The sequence $(\varepsilon^k)_{k \in {\Bbb N}}$ must then be
chosen as a sequence of type $\sigma$ (i.e. \hbox{$\sum_{k = 0}^{ k = +\infty} \varepsilon^k =
  +\infty$} and \hbox{$\sum_{k = 0}^{ k = +\infty} \left(\varepsilon^k\right)^2 < +\infty$}). The
convergence is necessarily slow.

Nevertheless, the non-differentiability of the dual function is not the main
difficulty. In this case, to find a price $p^{*}$ which maximizes the dual function $\Psi$ is not
enough~: for $p=p^{*}$, the primal minimization of the Lagrangian will not necessarily give a
solution of (\ref{pbtyp}) \cite{renIEEE}. In practise, a ``small'' variation of the ``prices'' leads
to a large variation of the primal variables. The primal variables ``switch'' from one value to
another  and never satisfy the coupling constraints. 
To our mind, these theoretical difficulties explain to a large extent the bad
reputation that these dual methods have in terms of convergence. 

However, Augmented Lagrangian can be used to reduce these
difficulties.
The Augmented Lagrangian $L_{c}$ related to the problem (\ref{pbtyp}) is defined over ${\cal U} \times {\cal C}$ as
follows~:
\begin{eqnarray}
  \forall (u,p) \, : \,  L_{c}(u,p) \eqdef L(u,p) + { c
  \over 2 } \Vert \Theta(u) \Vert^{2} \label{lagaug} 
\end{eqnarray} 
In the convex case, the	saddle-points of this Lagrangian are the same
as those of $L$ \cite{cohenzhu}. Then, the dual function
$\Psi_{c}$ related to $L_{c}$ is differentiable \cite{dual}.
Furthermore, solving $\min_{u \in U^{\ad}} L_{c}(u,p^{*})$ where $p^{*}$ is a maximum of the dual function necessarily yields a solution of
(\ref{pbtyp}). 

At first sight, this Augmented Lagrangian technique has a major
drawback with regards to decomposition : it introduces non-separable
terms ${ c \over 2 } \Vert \Theta(u) \Vert^{2}$. 
But, this difficulty can be overcome by linearizing the non-separable terms at each iteration~\cite{cohenzhu}. This
strategy leads to considering the following algorithm ({\bf algorithm A}) (iteration $k+1$):
\begin{itemize}
\item for all $i \in I$, $u^{k+1}_i$ is computed by solving~:
  $$\min_{u_i \in {\cal U}^{\ad}_i} \left( 
    J_i(u_i) +
    \left\langle \pi^k, \Theta_i(u_i)\right\rangle  +  { b \over 2 }
    \Vert u_{i} - u_{i}^{k} \Vert^{2}\right)\enspace,$$
\item $p^{k+1} = p^k + \varepsilon \Theta(u^{k+1})$,
\end{itemize}
with $b \in {\Bbb R}^{+*}$ and $ \pi^k= p^k + c \Theta(u^{k})$.

In the convex case, even if the cost function $J$ is not strictly convex, the convergence of this
algorithm towards a saddle point of $L$ has been proven provided that $0 < \varepsilon < 2c$ and $c
\tau^{2} < b$, where $\tau$ is the Lipschitz constant of $\Theta$~\cite{cohenzhu}.

This algorithm has been shown to be particularly efficient in dealing with classical hydro-thermal
generation scheduling problems~\cite{PSCC90},~\cite{renIEEE}. 
\subsubsection{Application to the hydro-scheduling}
\underline{\bf problem}
In the problem (\ref{pb}), two types of constraints have to be handled:
\begin{itemize}
\item state constraints ---~the volumes bounds: $V_i(t) \in [\Vmin{i},\Vmax{i}]$~---,
\item logical constraints concerning the controls which can be rather complex ---~$T \in \uad$~---~.
\end{itemize}
Suppose that dynamic constraints~(\ref{dynam}) do not have to be dealt with. It would not
be necessary to handle these two types of constraints simultaneously  and the decomposition of
(\ref{pb}) in simple subproblems would be allowed. 

This remark led us to
dualize dynamic constraints~(\ref{dynam}). 
Then, the application of the Algorithm~A to~(\ref{pb}) results, at the $k+1$
iteration, in the following steps\footnote{For the sake of completeness, it must be emphasized
  that this algorithm is not exactly the algorithm (A). The discharges which solve~(\ref{Alocturdua}) are
  used to define the cost function of the ``volume'' problems~(\ref{Alocexcdua}) which are solved at the same
  iteration. In practise, this sequential version turns out to be more efficient. }: 
\begin{itemize}
\item for all $i \in {\cal I}$, resolution of: 
  \begin{eqnarray}
    \lefteqn{ \displaystyle \displaystyle \min_{T_{i}([0,H[)}
    \dis \sum_{t \in [0,H[}\left\{ - c_i(T_i(t),t) + { b \over 2} (T_i(t) - T^k_i(t))^2 \enspace\enspace\right. } \nonumber\\*
    && \mbox{} + \left. \left(\pi^{k,k}_i(t) - \pi^{k,k}_{
       {\avi}}(t + \delta_{{\avi}\leftarrow i}) \right) T_i(t)\right\}\label{Alocturdua}\\*
    &&\mbox{s.t.\hspace{0.05cm}:}\enspace\enspace \turbi{T}\nonumber
  \end{eqnarray}
  yields ${T}^{k+1}_i([0,H[)$,
\item for all $i \in {\cal I}$ and $t \in [0,H[$, resolution of:
  \begin{eqnarray}
    \lefteqn{ \displaystyle  \dis \min_{{D}_{i}(t),{V}_{i}(t)} 
    \dis \left\{ -  \CfinI(V_i(t),t) \right.}\label{Alocexcdua}\\*
    && \mbox{}+ \dis \left(\pi^{k,k+1}_i(t) - \pi_{\avi}^{k,k+1}(t+ \delta_{{\avi} \leftarrow i}) \right) D_i(t) \nonumber\\*
    &&  \mbox{}+ \dis {b \over 2}  ({D}_i(t) - {D}^{k}_i(t))^{2}  
       +  {b \over 2}   ({V}_i(t) - {V}^{k}_i(t))^2 \nonumber\\*
    && \left.\mbox{}+ \dis \left(\pi^{k,k+1}_{i}(t-1) - \pi^{k,k+1}_{i}(t) \right)  V_i(t)\right\}\nonumber\\*
    &&\mbox{subject to:} \nonumber\\*
    &&V_i(t) \in [\Vmin{i},\Vmax{i}]\,,\quad 
       D_i(t) \in [0,\Dmax{i}]\, ,\nonumber\\*
    && \exclusion{V}{D} \,, \nonumber
  \end{eqnarray}
  yields  $\dis ( {V}^{k+1}_i(t), {D}^{k+1}_i(t))$,
\item for all  $i \in {\cal I}$ and for all $t \in [0,H[$ the dual variables $p_{i}(t)$ are updated
  as follows: \[
    p^{k+1}_{i}(t) =  p^k_{i}(t) + 
    c H_{i}^{k+1,k+1}(t)
  \]
\end{itemize}
where:
\begin{itemize}
\item $c>0$ and  $b>0$,
\item $\touti\,\;\forall t \in [0,H[~:$
  \begin{eqnarray}
    \lefteqn{ \displaystyle  H_{i}^{k_1,k_2}(t)
    \eqdef  V^{k_1}_i(t+1) - V^{k_1}_i(t) +
    D^{k_1}_i(t) } \nonumber\\*&&
                                  \mbox{} + T^{k_2}_i(t) - A^*_i(t) - \sum_{j \in \Gamma^T_i}
                                  [D^{k_1}_j+
                                  T^{k_2}_j](t_\jversi)  \;,\nonumber\\*
    \lefteqn{ \displaystyle  \pi^{k_1,k_2}_i(t) \eqdef  p^{k_1}_i(t) + c H_{i}^{k_1,k_2}(t)\enspace,}\nonumber
  \end{eqnarray}
\item $\forall t \in [0,H[\,:\, \CfinI(V_i(t),t) = 0 $.
\end{itemize} 
At each iteration, a dynamic subproblem related to each plant (\ref{Alocturdua}) is
solved. This subproblem handles mixed integer constraints concerning the plant discharge but no state
constraints. It is solved by dynamic programming. 

The problems (\ref{Alocexcdua}) are very
small and do not present any difficulties: only two real variables are optimized. 
Therefore, the subproblem resolutions that this algorithm requires, turn out to be
quite simple. 

Nevertheless, this is explained by the dualization of the most important
constraints of the problem (\ref{pb}): the dynamic constraints. It may seem dubious,
considering the mixed-integer constraints which have to be handled, that this dual method should
achieve a feasible solution of (\ref{pb}).

To explain the numerical result which will be outlined further, it may be emphasized that
this algorithm has been implemented in the following way:
\begin{itemize}
\item First of all, the ``convexified'' problem is solved using the
  price decomposition algorithm we have already described. In this case, the convexity assumptions
  being met, convergence is theoretically ensured and is obtained in practise. This first step
  yields a very good initial value of the dual variables $p_{i}(t)$.  
\item Then, every hundred iterations, until a feasible solution is
  reached, parameters are modified as follows:
  \begin{itemize}
  \item  the minimal bounds on the volumes $V_i^{\min}$ are slightly increased, 
  \item  the value of parameter $c$ of the Augmented Lagrangian is
    multiplied (by $3$). 
  \end{itemize}
  Even if it has not been explained theoretically, this progressive
  increasing of parameter $c$ turns out to be a very efficient method for obtaining feasible solutions.
\end{itemize}
\subsection{Interaction Prediction Principle}
The second decomposition technique we consider lies on a simultaneous partitioning 
of variables and constraints.
Every subproblem updates a set of variables handling a part of the constraints. Prices
remunerate the sharing in the satisfaction of the constraints which are not coped with. 
Hence, this approach may be considered as mixing the 
price and resources decomposition techniques.
\subsubsection{Theoretical Background}
\subparagraph{Takahara algorithm:}
Consider problem (\ref{pbtyp}). Suppose that  $\Theta= \prod_{i=1}^n \Theta_i $
where: for all $i \in I$ : $\Theta_i: U \rightarrow {\cal C}_i$ and ${\cal C}= \prod_{i=1}^n {\cal
  C}_i$.

At the iteration $k$, Takahara algorithm \cite{takahara}, \cite{takaconv} substitutes to (\ref{pb}) a sequence of
subproblems~(\ref{pbloctak}): 
\begin{equation}
  \begin{array}{l} \displaystyle \min_{u_i \in {\cal U}^{\ad}_i}
    J_i(u_i)+ \sum_{j\ne i} \psca{ p_j^k,
    \T{j}{u_i,u_{-i}^k} } \\
    {\ds \Theta_i(u_i,u_{-i}^k)=0 } 
  \end{array} 
  \label{pbloctak}	
\end{equation}
where $(u_i,u_{-i}^k)$ denotes the vector whose components are equal to those of $u^{k}$ except
$u_{i}$. 

Resolution of each (\ref{pbloctak})  yields a primal solution
$u_i^{k+1}$ and dual variables ($p_i^{k+1}$) related
to the local constraint $  \Theta_i(u_i,u_{-i}^k)=0 $. $(p_j^k)_{j\ne
  i}$ denote the dual variables that have been ``predicted'' by the
other subproblems at step $k$. They are used to remunerate the
participation of problem $i$ to the other constraints.
\subparagraph{A primal-dual relaxation algorithm:} To explain the
nature of this algorithm,
let us assume that $J$ is differentiable and ${\cal U}^{\ad}={\cal
  U}$. 
Then, Kuhn and Tucker necessary optimality conditions related to (\ref{pb}) may be written as follows:
\begin{equation}
  \left\{
    \begin{array}{l}
      \forall i \in I \;: \; J^{'}_i(u_i) + \left(\Theta'_{u_i}\right)^* p = 0 \enspace,\\
      \T{}{u}=0\enspace.
    \end{array}
  \right.\label{optimal}
\end{equation}
Furthermore, if $(u^{k+1}_i,p^{k+1}_i)$ is a (primal-dual) solution of (\ref{pbloctak}) then:
\[
  \left\{
    \begin{array}{l}
      J^{'}_i(u_i^{k+1}) + \sum_{j\ne i} \left(\Theta'_{j,u_i}\right)^* p^k_j +\left(\Theta'_{i,u_i}\right)^* p^{k+1}_i = 0 \enspace,\\
      \Theta_i{(u_i^{k+1},u_{-i}^k)}=0\enspace.
    \end{array}
  \right.
\]
Consequently, the Takahara algorithm appears to be a primal-dual relaxation algorithm applied to the resolution of (\ref{optimal}).
\subparagraph{Find a saddle-point of the Augmented Lagrangian:} 
The algorithm we apply to the hydro-problem is built up in this way. 
However, this primal-dual relaxation framework is not used to find a saddle-point of $L$ but 
a saddle-point of the Augmented Lagrangian $L_c$ (\ref{lagaug}). 

With the notation introduced above, at iteration $k$, it leads to solving the following 
subproblems:
\begin{eqnarray} 
  \begin{array}{l} \displaystyle \min_{u_i \in \U2ad_i} \left( \begin{array}{c}
                                                                 J_i(u_i)+ \sum_{j\ne i} \psca{ p_j^k,
                                                                 \T{j}{u_i,u_{-i}^k} } \\*
                                                                 + {c \over 2} \Vert \T{}{u_i,u_{-i}^k} \Vert^2
                                                               \end{array} \right)
    \\
    {\ds \Theta_i(u_i,u_{-i}^k)=0 } 
  \end{array} 
  \label{ploctakaug}	
\end{eqnarray}
Although we will not present hereafter a comparative test on this purpose, it may be pointed
out that, applied to the hydro-problem, (\ref{ploctakaug}) turns out to be more efficient than
(\ref{pbloctak}).
\subsubsection{Application to the hydro-problem} 
In order to apply this algorithm, we first reformulate
the hydro-problem (\ref{pbtyp}). Variables $A_i^+(t)$ representing the global inflows of
each reservoir are introduced:  

\begin{equation} A_i^+(t) -
  A_i^*(t) - \sum_{j \in \Gamma^T_i} [D_j+T_j](\tt{i}{j})=0
  \enspace. \label{inflows}
\end{equation}

With this definition, the constraint
\begin{eqnarray} 
  \displaystyle  V_i(t+1) \,  = \,  V_i(t) + A_i^+(t) - D_i(t) - T_i(t) \enspace,
  \label{enplusdyn} 
\end{eqnarray}
appears to be equivalent to the dynamic constraints~(\ref{dynam}).

To split (\ref{pb}) into a sequence of subproblems, each related to a plant, algorithm
(\ref{ploctakaug}) is then applied in the following
way: \begin{itemize} 
\item vector $\left(A_i^+(t), T_{i}(t), D_{i}(t) \right)_{t \in [0,H[}$ is $u_{i}$,
\item constraints (\ref{inflows}) are considered as being the coupling
  constraints $\Theta$, 
\item constraints (\ref{enplusdyn}), (\ref{exclus}), (\ref{borne}) 
  and (\ref{turbi}) define the domain
  ${\cal U}^{\ad}$ of (\ref{pbtyp}),
\item $J_i(u_{i})= \displaystyle \critere$.
\end{itemize}
With these choices, the algorithm (\ref{ploctakaug}) leads to solving,
at iteration $k+1$, the following subproblems  ({\bf Algorithm B}):
\begin{eqnarray}
  \displaystyle \lefteqn{ \displaystyle \min_{ T_i\in \uad_{i}} \sument \left\{
  \criteredans  \right.} \nonumber\\*&&
                                        \mbox{} - p_{\avi}^k(t)\,[ T_i + D_i](\tt{\avi}{i}) \nonumber\\*&&
                                                                                                           \left. \mbox{}	+ { c \over 2 } \norm{ DT_{i}^{k}(t) -
                                                                                                           [D_i+ T_i](\tt{\avi}{i})}^{2} \right\} \nonumber \\*
  \lefteqn{ \displaystyle  \mbox{s.t.: } 
  V_i(t+1) \,  = \,  V_i(t)  - D_i(t) - T_i(t) + A^*_i(t) } \nonumber
  \\*&&
        \mbox{} + \sum_{j \in \Gamma^T_{i}} \left[
        D^{k}_j+T^{k}_j\right](t_{i \leftarrow j})  \nonumber 
  \\*&&
        \mbox{\hspace{0.4cm}(\ref{exclus}), (\ref{borne}) and (\ref{turbi}) \enspace.}\nonumber
\end{eqnarray}
\begin{eqnarray}
  \mbox{with:}&&\quad DT_{i}^{k}(t)\, =\,
                 A_{\avi}^+(t)^{k} -
                 A_{\avi}^*(t)^{k} \nonumber \\* &&\displaystyle
                                                    \qquad - \sum_{j \in \Gamma^T_{\avi}-\{i\}}
                                                    [D^{k}_j+T^{k}_j](\tt{\avi}{j})
                                                    \enspace.\nonumber
\end{eqnarray} 
A sequential version of this algorithm has been implemented. Following
the course of the river, each subproblem is solved taking into account the results of the current
iteration (for the upstream informations) and of the preceding iteration (for the downstream). 
In this context, $DT_{i}^k(t)$ is necessarily equal to $[D^{k}_i + T^{k}_i](\tt{\avi}{i})$ and at iteration
$k+1$, the subproblem related to the plant $i$ may be written as follows:
\begin{eqnarray}
  \displaystyle \lefteqn{ \displaystyle  \min_{ T_i \in \uad_{i}} \sument \left\{
  \criteredans \right. } \label{loctakhyd}\\*&&
                                                \mbox{} - p_{\avi}^k(t)[ T_i + D_i](\tt{\avi}{i}) \nonumber \\*&&
                                                                                                                  \left. \mbox{} + { c \over 2 }{\Vert [D_i+ T_i]^k(\tt{\avi}{i}) - [D_i+ T_i]
                                                                                                                  (\tt{\avi}{i})\Vert}^{2} \right\}\nonumber \\*
  \lefteqn{ \mbox{s.t.: } 
  V_i(t+1) \,  = \,  V_i(t)  - D_i(t) - T_i(t)  + A^*_i(t)} \nonumber\\*&&
                                                                           \mbox{}\displaystyle\hspace{0.5cm} + \sum_{j \in \Gamma^T_{i}}
                                                                           [D^{k+1}_j+T^{k+1}_j](\tt{i}{j})\nonumber\\*&& 
                                                                                                                          \mbox{\hspace{0.4cm}(\ref{exclus}), (\ref{borne}) and (\ref{turbi})}\enspace. \nonumber
\end{eqnarray}
Therefore, handling its generating constraints and taking into
account the discharge the upstream plants have computed, each plant
optimizes its schedule.
The dual variable $p_{\avi}^k$ may be understood as being the price the downstream reservoir
``would pay'' its water inflows. Quadratic terms introduced by  the Augmented Lagrangian appear to be a
type of ``brake'' avoiding the oscillations of the algorithm.

It may also be noticed the subproblems (\ref{loctakhyd}) turn out to
have exactly the same structure 
as the local problems of
the heuristic process currently in use (See Today's methods). Therefore, from a practical point of
view, Algorithm~B appears to be an extension of this heuristic method. Moreover, providing there is no pumping unit, this sequential version
{\it generally}\footnote{generally but not always} yields a feasible solution at the first iteration.

\section{Numerical Tests}
\subsection{A hydraulic valleys sample}
A sample of hydraulic valleys has been chosen so as to point out, as well as
possible, the main resolution difficulties. 

\begin{figure}[hbtp]
  \begin{center}
    \includegraphics{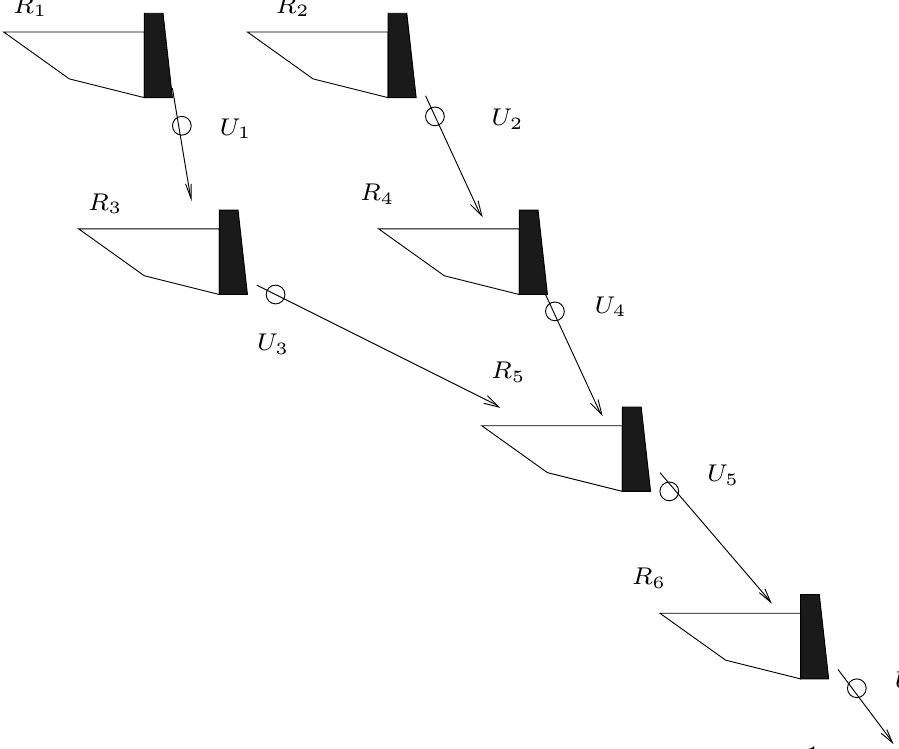}
  \end{center}
  \caption{\label{vallee}Hydraulic test valley}
\end{figure}

{\bf Size and topology:} The hydraulic valley is illustrated in
Figure~\ref{vallee}. It contains six reservoirs and to each reservoir $R_i$ is related a plant $U_i$.

{\bf Volume bounds:} The upstream reservoirs 1 and 2 are supposed to have a large storage
capacity: on a daily scale, no volume constraints have to be coped with. 
The other reservoirs are characterized by the ratio of their storage and hourly-discharge
capacities. Three sets of storage/discharge ratios are considered:
\begin{itemize}
\item [{\bf V1}] $1\mbox{ h }30$ for reservoirs $2$ and $4$, $5\mbox{h}$ for reservoirs $5$ and
  $6$,
\item [{\bf V2}] $0\mbox{ h }30$ for reservoirs $2$ and $4$, $3\mbox{h}$ for reservoirs $5$ and
  $6$,
\item [{\bf V3}] $5\mbox{ h }$ for reservoirs $2$ and $4$, $15\mbox{h}$ for reservoirs $5$ and
  $6$.
\end{itemize}
There are no natural inflows. 

\begin{figure}[hbtp]
  \begin{center}
    \includegraphics[width=\textwidth]{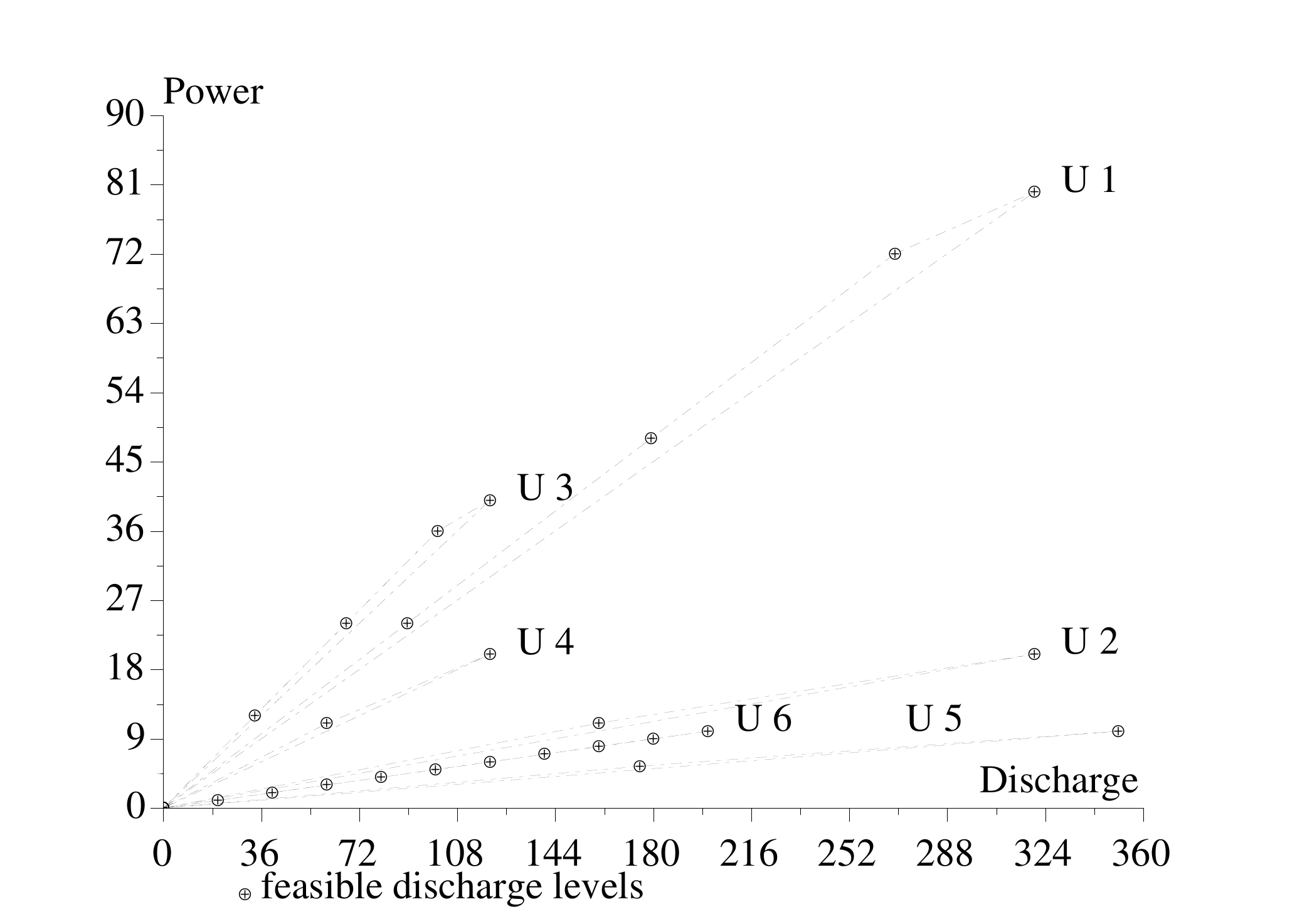}
  \end{center}
  \caption{\label{fig:usine}The generation is a piecewise linear function of the discharge}
\end{figure}

{\bf Cost function:} The generation ``revenue'' $c_i(T_i(t),t)$ is
assumed to be a linear function of the discharge (i.e. $c_i(T_i(t),t)=
p^{\gen}(t) T_{i}(t)$). Three sequences $(p^{\gen}(t))_{t \in [0,H[}$ are considered (in Francs per MWh):
\begin{itemize}
\item [{\bf P1}] The first ranges from $99$ to $101$. It ``switches'' from one value to another every
  $4$ hours. This price vector enables the numerical accuracy of the algorithms to be tested.  
\item [{\bf P2}] The second remains at $100$ over the whole period except four hours during which it
  rises up to $500$. Such a choice enables the spinning reserve over a four hour period to be computed.
  The ability to optimize {\it feasible} controls in
  real-time and in case of emergency is also measured in this way. From a numerical point of view, it
  is in this case that constraints (\ref{exclus}) and (\ref{borne}) are actually active. 
\item [{\bf P3}] Every four hours, the third switches between $80$ and $120$. 
\end{itemize}
{\bf Discharge Constraints: }For each plant $U_i$, the discharge belongs to a discrete set of values (See
Figure~\ref{fig:usine}). 

The minimum delay between two discharge variations may be $0$ ({\bf D1}), $1$ ({\bf D2}) or $4$
({\bf D3}) hours. 

The ({\bf D4}), ({\bf D5}) and ({\bf D6}) are deduced from ({\bf D1}), ({\bf D2}) or
({\bf D3}) by supposing the plants $3$ and $4$ are out-of-order (have no discharge capacity).

Crossing these factors, it is a $56$ valleys sample which is built up.
It may be noticed the
non-continuities in the discharge domain are actually sizeable. For the {\bf V2}
storage/discharge ratio, the choice of a discharge level rather than another modify the hourly discharge by about half the storage capacity~!
\subsection{Numerical results} The two decomposition algorithms and the process
currently in use at EDF have been compared over this sample of valleys. 

The ``scores'' are computed in the following way. For each tested
algorithm and each sample of valley, the maximum gain $\gain$ obtained
by a feasible solution along the iterations is recorded.
If no feasible solution is found, this gain is considered to be zero: because there are no natural
inflows, a zero discharge solution is possible and its gain is zero. 

Then, relative gains (or ``scores'') are computed by dividing these gains by those
of the best solutions achieved by one of the three processes. 

\begin{table}[hbtp]
  \center{
    \begin{tabular}{|l||p{1.0cm}|p{1.0cm}|p{1.0cm}|p{1.2cm}|}\hline
      & Alg. & Alg. & Alg. & Alg.  \\ 
      &{\bf A}& {\bf B}&{\bf C}& {\bf L} \\ 
      \hline\hline
      {\bf P1} (m.s.) & 83.2 \% 	& 83.8 \% 	&  62.3 \% 	&  186 \% \\ \hline
      {\bf P1} (a.f.) & 100 \% 	& 100\% 	&  72  \% 	& 0 \%\\ \hline\hline
      {\bf P2} (m.s.) & 99.4 \% 	& 99.9 \% 	&  99.7  \% 	& 100 \% \\ \hline
      {\bf P2} (a.f.) & 100  \% 	& 100\% 	& 100 \% 	& 0 \%\\ \hline\hline
      {\bf P3} (m.s.) & 97.1\% 	& 90.4  \% 	& 90.8 \% 	& 109.2 \% \\ \hline
      {\bf P3} (a.f.) & 100\% 	& 100 \% 	& 100 \% 	& 0 \%\\ \hline\hline\hline
      {\bf M} (m.s.) & 93.2 \% 	& 91.4  \% 	& 84.3 \% 	& 131.8 \% \\ \hline
      {\bf M}  (a.f.) & 100\% 	& 100 \% 	& 91  \% 	& 0 \%\\ \hline
    \end{tabular}}
  \caption{Mean scores (m.s.) and average feasibility (a.f.) for each 
    ``prices system'' ({\bf P1}, {\bf P2}, {\bf P3}) and for the whole sample {\bf M}.
    \label{lesmoyennes}}
\end{table}

Table~\ref{lesmoyennes} gives mean scores (m.s.) and feasibility
average rates (f.a.) for the three types of
``generation prices'' which have been considered ({\bf P1}, {\bf P2}, {\bf P3}) and for the whole
sample ({\bf M}). In every case, the two decomposition methods yield
feasible solutions. Considering how significant the non-continuities
are it is a remarkable achievement. 

In spite of the efficiency of the heuristic process which is currently in use at EDF ({\bf Algorithm C}) --- in more than
$90$
percent of the (difficult) cases we have selected, feasible solutions are reached --- these methods
thus represent a real improvement. 

Furthermore, if constraints~(\ref{turbi}) and~(\ref{exclus}) are not handled, problem
(\ref{pb}) is convex. It can be solved by linear programming. The mean ``score'' of the non feasible
solutions obtained in this way are also indicated in
Table~\ref{lesmoyennes} ({\bf Algorithm L}). One may consider that
about twenty percents of the current cost of the
nonconvexities in the (\ref{pb}) modelling are saved thanks to these decomposition methods.

\paragraph{What CPU time~?}
For the valley of six reservoirs we consider, the heuristic resolution today in use at EDF takes about
$10$ seconds on a SUN $4/40$\footnote{All the CPU time here mentioned have been measured on this
  computer}. One half of this time is dedicated to the linear optimization the other is used by the six
dynamic programming resolution which are necessary to find a solution. 

For the prediction strategy, the results presented above correspond to $25$ iterations. Each of these
having the same complexity as the heuristic research of a feasible solution, the CPU time
required is more or less $2$ minutes. 

Our implementation of the price decomposition method uses about $1200$ iterations. This number may
seem important. Nevertheless, the subproblems are particularly simple and, for the six reservoirs
valleys, CPU time does not exceed $2$ minutes. 

If these times are not huge, they multiply by ten CPU times of the current process. Therefore, work is
currently undertaken to reduce these CPU times. To our mind, on average, they should be divided by
about $5$ in the final implementation by:
\begin{itemize}
\item improving the software design,
\item avoiding to solve each subproblem at each iteration.
\end{itemize}
\paragraph{What is the best method~?}
We have already noticed the average ``score'' ({\bf M}) of the price decomposition is $93.2\%$, the
prediction one being $91.4\%$. Should the prediction method be rejected~? We have not made such a
choice for two reasons:
\begin{itemize}
\item contrary to the price decomposition strategy, the
  relaxation  algorithm generally yields, from the first
  iterations, a feasible solution, 
\item if, on average, price decomposition method 
  reaches the best solutions, it does not in every case.
  Furthermore,
  if one choose the best of the
  two solutions these methods yield, it would not be $93.2\%$ or $91.4\%$
  but a score of $100\%$ which would be reached. In fact, the tools we
  are currently developing,
  on the basis of these first tests, will try to take advantage of each
  of these methods.
  By experimentations, we aim at establishing rules which,
  after considering the characteristics of the valley, choose the best of 	the two algorithms. 
\end{itemize}
\section{Conclusion} Over the $56$ numerical tests which have been
carried out, in spite of the mixed-integer constraints which are handled, the decomposition methods
considered allow a feasible solution to be systematically found. Moreover, compared to the two-step
process currently in use at EDF, those methods yield sizeable savings. 

For these reasons, these methods will be used to design new regional level software at
Electricit\'{e} de France.

Furthermore, by proving the
robustness of these decomposition approaches, these tests open up new fields of
research. A global optimization of generation schedules of several hydraulic valleys, handling
coupling constraints (demand constraints), could be achieved in this
way.

\end{document}